\documentclass{article}
\usepackage{amsmath}
\usepackage{amsthm}
\usepackage{amssymb}
\usepackage{graphics}
\usepackage{graphicx}
\usepackage{enumerate}
\usepackage{vaucanson-g}

\newtheorem{thrm}{Theorem}[section]
\newtheorem{lmm}[thrm]{Lemma}

\theoremstyle{definition}
\newtheorem{rmrk}[thrm]{Remark}

\def\N{\mathbb N}

\def\C{\mathcal C}
\def\pf{\begin{proof}}
\def\pfk{\end{proof}}

\numberwithin{equation}{section}

\def\N{\mathbb N}

\def\C{\mathcal C}

\newtheorem{obs}[thrm]{Observation}
\def\pf{\begin{proof}}
\def\pfk{\end{proof}}
\begin{document}
\title{Palindromes in infinite ternary words}

\author{L\!'ubom\'ira Balkov\'a $^{(1)}$
\and Edita Pelantov\'a $^{(1)}$
\and \v St\v ep\'an Starosta $^{(1)}$ $^{(2)}$\footnote{{\em E-mail addresses}: l.balkova@centrum.cz, edita.pelantova@fjfi.cvut.cz, starosta@iml.univ-mrs.fr}}

\maketitle
\thispagestyle{empty}

\begin{center}
(1) Doppler Institute for Mathematical Physics and Applied Mathematics and\\
Department of Mathematics, FNSPE, Czech Technical University, \\
Trojanova 13, 120 00
Praha 2, Czech Republic\\
(2) Institut de Math\'ematiques de Luminy, Campus de Luminy, \\ Case 907, 13288 MARSEILLE Cedex 9
\end{center}


\noindent {\em Date}: January 27, 2009

\begin{abstract}
We study infinite words $\mathbf{u}$ over an alphabet $\mathcal{A}$
satisfying the property
 $$\mathcal{P} : \qquad   \mathcal{P}(n)+
\mathcal{P}(n+1) = 1+ \#\mathcal{A} \quad \hbox{for any} \ n \in
\mathbb{N},  $$  where $\mathcal{P}(n)$  denotes the number of
palindromic factors of length $n$ occurring in the language of~$\mathbf{u}$.
We study also infinite words satisfying a~stronger property
$$\mathcal{PE}: \quad \hbox{every palindrome of $\mathbf{u}$ has
exactly one palindromic extension in $\mathbf{u}$}\,.$$  For
binary words, the properties $\mathcal{P}$ and $\mathcal{PE}$
coincide and these properties characterize Sturmian words, i.e.,
words with the complexity $\mathcal{C}(n)=n+1$ for any $n
\in \mathbb{N}$. In this paper, we focus on ternary infinite words
with the language closed under reversal. For such words $\mathbf{u}$,
we prove that if $\mathcal{C}(n)=2n+1$ for any $n \in \mathbb{N}$
then $\mathbf{u}$ satisfies the property $\mathcal{P}$ and
moreover $\mathbf{u}$ is rich in palindromes.
 Also a~sufficient condition for the property $\mathcal{PE}$ is given.
We construct a~word demonstrating that $\mathcal{P}$ on a~ternary
alphabet does not imply $\mathcal{PE}$.
\end{abstract}

\section{Introduction}
Sturmian words are the most intensively studied infinite words
since their appearance in 1940. They were introduced by Morse and
Hedlund~\cite{HeMo} as aperiodic words with the minimal possible
complexity, i.e., with the complexity
$\mathcal{C}(n)=n+1$ for any $n\in \mathbb{N}$. The complexity is the function $\mathcal{C}: \mathbb{N}\mapsto
\mathbb{N}$ defined by
$$ \mathcal{C}(n)=\hbox{number of  factors of length $n$ occurring in} \ \mathbf{u}
\,.$$ The set of all factors occurring in~$\mathbf{u}$ is called the
{\em language} of $\mathbf{u}$ and denoted throughout this paper by
$\mathcal{L}(\mathbf{u})$. There exist many equivalent definitions
of Sturmian words. Already in~\cite{HeMo}, Sturmian
words are characterized by their balance property. In the center
of our attention will be another characteristics  of Sturmian
words, recently proved in~\cite{DrPi}. This characteristics uses the
palindromic complexity of~$\mathbf{u}$, which is the
function $\mathcal{P}: \mathbb{N}\mapsto \mathbb{N}$ defined by
$$ \mathcal{P}(n)=\hbox{number of palindromic factors of length $n$ occurring in} \ \mathbf{u}
\,.$$ Droubay and Pirillo proved that an infinite word $ \mathbf{u}$ is Sturmian if and only if its
palindromic complexity is
$$  \mathcal{P}(n)  = \left\{\begin{array}{cl}
1 & \hbox{if}\ \  n \ \  \hbox{is even},\\
2 & \hbox{if}\  \ n \ \  \hbox{is odd}.\\
\end{array}\right.
$$
Since the empty word is the only palindrome of length $0$ and the letters
of the alphabet~$\mathcal{A}$ are the only palindromes of length $1$ in~$\mathbf{u}$,
the previous property can be rewritten in a~compact form for binary infinite words as
$$\qquad \mathcal{P}(n)+ \mathcal{P}(n+1)=3\quad \hbox{for any} \ n \in
\mathbb{N}.$$

Being inspired by Sturmian words, we generalize the previous property
for infinite words over any alphabet $\mathcal A$ as
$$\mathcal{P}: \qquad \mathcal{P}(n)+ \mathcal{P}(n+1) = 1+
\#\mathcal{A} \quad \hbox{for any} \ n \in
\mathbb{N}.~~~~~~~~~~~~~~~~~~~~~ $$
It is again readily seen that the property $\mathcal{P}$ is equivalent with the property
$$  \mathcal{P}(n)  = \left\{\begin{array}{cl}
1 & \hbox{if}\ \  n \ \  \hbox{is even},\\
\#\mathcal{A} & \hbox{if}\  \ n \ \  \hbox{is odd}.\\
\end{array}\right.
$$

Examples of infinite words over multilateral alphabets satisfying the property $\mathcal P$
are Arnoux-Rauzy words and nondegenerate words coding the $r$-interval
exchange transformation with the permutation $\pi =(r,r-1,r-2,\ldots,2,1)$.

When studying in details the proof of Droubay and Pirillo,
we learn that a~binary word $\mathbf{u}$ is Sturmian if and only if $\mathbf{u}$
satisfies the following condition
$$  \mathcal{PE}: \qquad  \hbox{any palindromic factor of} \ \
{\mathbf{u}}\ \  \hbox{has a~unique palindromic extension in }
{\mathbf{u}}\,.
$$
In other words, for any palindrome $p \in \mathcal{L}(\mathbf{u})$
there exists a~unique letter $a\in \mathcal{A}$ such that $apa \in
\mathcal{L}(\mathbf{u})$. In fact, our two examples of words with
the property $\mathcal P$ -  namely Arnoux-Rauzy words and words
coding interval exchange  -  have even the property $
\mathcal{PE}$.

Infinite words over a~multilateral alphabet satisfying the
property $\mathcal{P}$ or $ \mathcal{PE}$ may be understood as one
of the possible generalizations of Sturmian words. It is evident
that $\mathcal{PE}$ implies $\mathcal{P}$. The inverse implication
holds over a~binary alphabet, but it need not hold in general.
The validity of $\mathcal{P}$ or $ \mathcal{PE}$ guarantees that the
language $\mathcal{L}(\mathbf{u})$ contains infinitely many
distinct palindromic factors. Such a~language need not be closed
under reversal.
Nevertheless in the sequel, we concentrate on the study of ternary
words whose language is closed under reversal. It is readily seen that
such words are recurrent and their Rauzy graphs have
a~non-trivial automorphism that will serve as a~powerful tool in
our consideration.

We will prove the following two theorems:

\begin{thrm}\label{kdyP} An infinite ternary word whose language is closed under
reversal has the property $\mathcal P$ if its complexity
satisfies $\mathcal{C}(n)=2n+1$.\end{thrm}

For the description of $\mathcal{PE}$, an important role is played by the notion of
a~left special factor: a~factor $w\in \mathcal{L}(\mathbf{u})$
is called {\em left special} if there exist at
least two different letters $a, b$ such that both $aw\in
\mathcal{L}(\mathbf{u})$ and $bw\in \mathcal{L}(\mathbf{u})$.
A~left special factor $w$ is called {\em maximal} if for
any letter $c \in \mathcal{A}$, the factor $wc$ is not left
special.

\begin{thrm}\label{kdyPE} An infinite ternary word $\mathbf{u}$ whose language is closed under
reversal has the property $\mathcal{PE}$ if its complexity
satisfies $\mathcal{C}(n)=2n+1$ and $\mathbf{u}$ has no maximal left
special factor. \end{thrm}

It is interesting to mention two corollaries of the previous theorems.
Vuillon~\cite{Vu} showed that a~binary infinite word is Sturmian if and only if each
of its factors has exactly two return words, i.e., Sturmian words
are precisely binary words satisfying the property
$$\mathcal{R}:  \quad \hbox{any factor of $\mathbf{u}$ has exactly $\#\mathcal{A}$ return
words}.
$$
In the paper~\cite{BaPeSt}, it is shown that a~ternary infinite
uniformly recurrent word $\mathbf{u}$ has the property
$\mathcal{R}$ if and only if its complexity satisfies
$\mathcal{C}(n)=2n+1$ and $\mathbf{u}$ has no maximal left special
factor. Consequently, for ternary infinite  words with the
language closed under reversal, $\mathcal{R}$ implies
$\mathcal{PE}$.

Theorem~\ref{kdyP} says that for infinite words whose language is
closed under reversal and whose complexity satisfies
$\mathcal{C}(n)=2n+1$, the following equation holds
\begin{equation}\label{amy}
 \mathcal{P}(n)+ \mathcal{P}(n+1) =  2+
\mathcal{C}(n+1)- \mathcal{C}(n)\,.
\end{equation}
Infinite words
fulfilling the above equation are in a~certain sense the richest
in palindromes, since according to~\cite{BaMaPe}, any infinite
word whose language is closed under reversal satisfies
$$ \mathcal{P}(n)+ \mathcal{P}(n+1) \leq 2+ \mathcal{C}(n+1)-
\mathcal{C}(n)\,.$$

Different descriptions of rich words defined by the equation
\eqref{amy} can be found in \cite{GlJuWiZa}. Infinite
ternary words with the language closed under reversal and the
complexity $\mathcal{C}(n)=2n+1$ form a~further class of rich
words.

In Section~\ref{prelim}, we recall basic notions from combinatorics on words. Section~\ref{proofs}
contains the proofs of Theorems~\ref{kdyP} and~\ref{kdyPE}. Section~\ref{counterex} provides two examples of words: the first one shows that the properties $\mathcal P$ and ${\mathcal{PE}}$ are not equivalent and the second one proves that the implications in Theorems~\ref{kdyP} and~\ref{kdyPE} cannot be reversed.

\section{Preliminaries}\label{prelim}

By $\mathcal{A}$ we denote a~finite set of symbols, usually called
{\em letters}; the set $\mathcal{A}$ is therefore called an {\em
alphabet}. A finite string $w=w_0w_1\ldots w_{n-1}$ of letters of
$\mathcal{A}$ is said to be a~{\em finite word}, its length is
denoted by $|w| = n$. Finite words over  $\mathcal{A}$ together
with the operation of concatenation and the empty word $\varepsilon$
as the neutral element form a~free monoid $\mathcal{A}^*$. The
assignment
$$w=w_0w_1\ldots w_{n-1} \quad \mapsto \quad \overline{w} =
w_{n-1}w_{n-2}\ldots w_{0}$$ is a~bijection on $\mathcal{A}^*$,
the word $\overline{w}$ is called the {\em reversal} or the {\em
mirror image} of $w$. A~word $w$ which coincides with its mirror
image is a~{\em palindrome}.

Under an {\em infinite word} $\mathbf u$ we understand an infinite string
$\mathbf{u}=u_0u_1u_2\ldots $ of letters from $\mathcal{A}$.
A~finite word $w$ is a~{\em factor} of a~word $v$ (finite or
infinite) if there exist words $w^{(1)}$ and $w^{(2)}$ such that
$v= w^{(1)}w w^{(2)}$. If $w^{(1)} = \varepsilon$, then $w$ is said
to be a~{\em prefix} of $v$, if $w^{(2)} = \varepsilon$, then $w$ is
a~{\em suffix} of~$v$. We say that a~prefix, a~suffix is {\em proper} if it is not equal
to the word itself.
The {\em language} $\mathcal{L}(\mathbf{u})$ of an infinite word
$\mathbf{u}$ is the set of all its factors. The factors of
$\mathbf{u}$ of length $n$ form the set denoted by
$\mathcal{L}_n(\mathbf{u})$. Using this notation, we may write
$\mathcal{L}(\mathbf{u})=\cup_{n\in
\mathbb{N}}\mathcal{L}_n(\mathbf{u})$. We say that the language
$\mathcal{L}(\mathbf{u})$ is {\em closed under reversal} if
$\mathcal{L}(\mathbf{u})$ contains with every factor $w$ also its
reversal $\overline{w}$.

For any factor $w\in \mathcal{L}(\mathbf{u})$, there exists an
index $i$ such that $w$ is a prefix  of  the infinite word
$u_iu_{i+1}u_{i+2} \ldots$. Such an index $i$ is called an {\em
occurrence} of $w$ in $\mathbf{u}$. If each factor of $\mathbf{u}$
has at least two occurrences in $\mathbf{u}$, the infinite word
$\mathbf{u}$ is said to be {\em recurrent}. It is easy to see that
if the language of $\mathbf{u}$ is closed under reversal, then
$\mathbf{u}$ is recurrent.

The {\em complexity} of an infinite word $\mathbf{u}$ is a
map $\mathcal{C}: \mathbb{N} \mapsto \mathbb{N}$, defined by
$\mathcal{C}(n)=\# \mathcal{L}_n(\mathbf{u})$. To
determine the increment of the complexity, one has to
count the possible {extensions} of factors of length $n$. A~{\em
right extension} of $w \in \mathcal{L}(\mathbf{u})$ is any letter
$a\in \mathcal{A}$ such that $w a\in \mathcal{L}(\mathbf{u})$.
The set of all right extensions of a~factor $w$ will be denoted by
${\rm Rext}(w)$. Of course, any factor of $\mathbf u$ has at least one right
extension. A~factor $w$ is called {\em right special} if $w$ has
at least two right extensions. Clearly, any suffix of a~right
special factor is right special as well. A~right special factor
$w$ which is not a~suffix of any longer right special factor is
called a~{\em maximal right special} factor.
Similarly, one can define a~{\em left extension}, a~{\em left
special} factor and ${\rm Lext}(w)$. We will deal only with
recurrent infinite words $\mathbf{u}$. In this case, any factor of
$\mathbf{u}$ has at least one left extension. If $a \in {\mathcal A}$ and $p$
is a~palindrome and $apa \in {\mathcal L}(\mathbf u)$, then $apa$ is said to be a~{\em palindromic extension} of $p$.
We say that $w$ is
a~{\em bispecial} factor if it is right and left special.
The role of bispecial factors for the computation of the
complexity can be nicely illustrated on Rauzy graphs.

Let $\mathbf{u}$ be an infinite word and $n\in\N$. The {\em Rauzy graph}
$\Gamma_n$ of $\mathbf u$ is a~directed graph whose set of vertices is
${\mathcal L}_n(\mathbf{u})$ and set of edges is ${\mathcal
L}_{n+1}(\mathbf{u})$. An edge $e\in{\mathcal L}_{n+1}(\mathbf u)$ starts
at the vertex $x$ and ends at the vertex $y$ if $x$ is a prefix
and $y$ is a suffix of $e$.
\begin{figure}[ht]
\begin{center}
\begin{picture}(220,30)
\put(50,20){\circle*{5}} \put(210,20){\circle*{5}}
\put(60,20){\vector(1,0){140}}
%\qbezier(50,20)(120,40)(190,20)
\put(10,8){$x=w_0w_1\cdots w_{n-1}$} \put(170,8){$y=w_1\cdots
w_{n-1}w_n$} \put(82,25){$e=w_0w_1\cdots w_{n-1}w_n$}
\end{picture}
\end{center}
\caption{Incidence relation between an edge and vertices in a
Rauzy graph.} \label{f}
\end{figure}
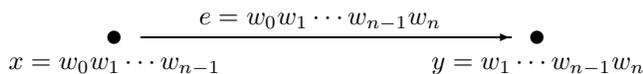
If the word $\mathbf{u}$ is  recurrent, the graph $\Gamma_n$ is
strongly connected for every $n\in\N$, i.e., there exists a~directed
path from every vertex $x$ to every vertex $y$ of the
graph.

The {\em outdegree} ({\em indegree}) of a~vertex $x\in{\mathcal
L}_n(\mathbf{u})$ is the number of edges which start (end) in $x$.
Obviously the outdegree of $x$ is equal to  $\#{\rm Rext}(x)$ and
 the indegree of $x$ is $\#{\rm Lext}(x)$.

The sum of outdegrees over all vertices is equal to the number of
edges in every directed graph. Similarly, it holds for indegrees.
In particular, for the Rauzy graph we have
$$
\sum_{x\in{\mathcal L}_{n}(\mathbf{u})}\#{\rm Rext}(x) \ =\
\mathcal{C}(n+1) \ = \sum_{x\in{\mathcal L}_{n}(\mathbf{u})}\#{\rm
Lext}(x)\,.
$$
The first difference of complexity $\Delta\C(n) = \mathcal{C}(n+1) - \mathcal{C}(n)$ is thus given by
\begin{equation}\label{eq:delta}
\Delta\C(n) = \sum_{x\in{\mathcal L}_{n}(\mathbf{u})}\bigl(\#{\rm
Rext}(x) - 1 \bigr) \ = \sum_{x\in{\mathcal
L}_{n}(\mathbf{u})}\bigl(\#{\rm Lext}(x) - 1 \bigr)\,.
\end{equation}
Let us restrict our consideration to recurrent words, then
a non-zero contribution to $\Delta\C(n)$ is given only
by those factors $x\in{\mathcal L}_n(\mathbf{u})$, for which
$\#{\rm Rext}(x)\geq 2$ or  $\#{\rm Lext}(x)\geq 2$, i.e., for
right or left special factors.  The relation~\eqref{eq:delta} can
be rewritten as
$$
\Delta\C(n) = \sum_{x\in{\mathcal L}_{n}(\mathbf{u}),\
\hbox{\scriptsize $x$ right
special}\hspace*{-1.5cm}}\hspace*{-0.3cm}\bigl(\#{\rm Rext}(x) - 1
\bigr) \  = \sum_{x\in{\mathcal L}_{n}(\mathbf{u}),\
\hbox{\scriptsize $x$ left
special}\hspace*{-1.5cm}}\hspace*{-0.3cm}\bigl(\#{\rm Lext}(x) - 1
\bigr) \ \,.
$$

If the language of the infinite word $\mathbf{u}$ is closed under
reversal, then  the operation that to every vertex $x$ of the
graph associates the vertex $\overline{x}$ and to every edge $e$
associates $\overline{e}$ maps the Rauzy graph $\Gamma_n$ onto
itself. In this case, we will draw the Rauzy graph $\Gamma_n$ axially symmetric
in the plane: the positions of vertices $x$ and $\overline x$ are
symmetrical with respect to an axis. Thus, $x$ is a~palindrome if and only if the vertex $x$ lies on the axis, and $e$ is a~palindrome of length $n+1$ if and only if the edge $e$ crosses the axis.

\section{Proof of Theorems~\ref{kdyP} and~\ref{kdyPE}}\label{proofs}
The proofs of Theorems~\ref{kdyP} and~\ref{kdyPE} will be a~consequence of the following three
lemmas that determine the number of palindromic extensions of palindromic factors with respect to the number of their
left extensions.

\begin{lmm}\label{Lemma1} Let  $\mathbf{u}$ be an
infinite word over an alphabet $\mathcal{A}$ whose language is
closed under reversal. If a~palindrome $p\in
\mathcal{L}(\mathbf{u})$ is not left special (and thus neither right special),
then $p$ has a~unique palindromic extension.
\end{lmm}

\pf Since $p\in \mathcal{L}(\mathbf{u})$ is not left special, there exists a~unique
$x \in \mathcal{A}$ such that $xp\in \mathcal{L}(\mathbf{u})$. By reversal closeness, $\mathcal{L}(\mathbf{u})$ contains also $px$. As $p$ has a~unique left extension $x$, the factor $px$ has $x$ as its unique left extension, too. Thus $xpx$ is the unique palindromic extension of $p$.
\pfk

\begin{lmm}\label{Lemma2} Let $\mathbf{u}$ be an infinite word over a~ternary
alphabet $\mathcal A$ with the complexity $\mathcal{C}(n) = 2n+1$ for
any $n \in \mathbb{N}$ and with the language ${\mathcal L}(\mathbf{u})$ closed under reversal.
If a~palindrome  $p \in \mathcal{L}(\mathbf{u})$ has $\#{\rm
Lext}(p)=3$, then $p$ has a~unique palindromic extension.
\end{lmm}

\pf  As $\Delta \mathcal{C}(n)=2$, the palindrome $p$ is the only left special factor
of length $n=|p|$, and by reversal closeness,  the only right special factor of length $n$, too.
\begin{enumerate}
\item  First, assume that there exists a~letter $x$ such
that ${\rm Lext}(px) = \mathcal{A}$. It means that $xpx$ is a~factor of $\mathbf{u}$, hence $xpx$ is a~palindromic extension of $p$. If there exists another palindromic extension of $p$,
i.e., $ypy \in {\mathcal L}(\mathbf{u})$ for $y\neq x$, then since $y\in {\rm
Lext}(px)$, it follows that $ypx$ and $xpy$ belong to $\mathcal{L}(\mathbf{u})$. Therefore $x,y \in {\rm Lext}(py)$,
which implies
 $$\Delta \mathcal{C}(n+1)\  \geq\  \#{\rm Lext}(px) - 1
+ \#{\rm Lext}(py)-1\  \geq\  3$$ - a contradiction.

 \item Second, suppose that for every letter $x$, it holds that
 ${\rm Lext}(px) \neq  \mathcal{A}$.  Let us recall that if
$w$ is a~left special factor of length $n+1$, then its prefix of length $n$ is necessarily
left special, too. As a~consequence, together with the fact that $\Delta C(n+1)=2$, there exist two left special factors $px$ and $py$ in
 $\mathcal{L}_{n+1}(\mathbf{u})$ for $x\neq y$ with $\#{\rm
Lext}(px)= \#{\rm Lext}(py)= 2$. Denote ${\rm Lext}(px)=\{a,b\}$
and ${\rm Lext}(py)=\{A,B\}$. Since our alphabet is
ternary, we may assume WLOG that $a=A$. By reversal closeness, it follows that
$xpa$ and $ypa$ belong to the language, and therefore, the factor $pa$ is left special as
well. WLOG $a=x$ and $b=y$. Denote by $c$ the third letter of $\mathcal A$.
Since $c \in {\rm Rext}(p)$ and by recurrence of $\mathbf{u}$, there exists a~letter $C$ such that $Cpc \in {\mathcal L}(\mathbf{u})$. However, since $pc$ is not left special, this $C$ is unique.

-- If $C=a$, then ${\rm Lext}(pa)=\mathcal{A}$ - a
contradiction.

-- If $C=b$, then necessarily $B=c$ and $apa$ is the unique
palindromic extension of $p$, as claimed.

-- If $C=c$, then $B=b$. The Rauzy graph $\Gamma_n$ has a~unique vertex of indegree $>1$ (see Figure~\ref{RauzyPicture}, where the straight lines denote edges and the zig zag lines denote paths) -- the bispecial factor $p$. Consequently, the vertex $p$ is the unique common vertex of three cycles. Since ${\rm Lext}(cp)=\{c\}$, after coming to the vertex $p$ using the edge $cp$, we cannot leave $p$ but using the edge $pc$. Hence,
we move eventually in a~unique cycle and the word $\mathbf{u}$ is thus eventually periodic - a~contradiction with the complexity.

\begin{figure}[!h]
\begin{center}
\MediumPicture\VCDraw{%
\begin{VCPicture}{(-5,-5)(5,5)}
\ChgStateLineWidth{0.5}\StateVar[p]{(0, 0)}{p} \VSState{(2, 4)}{0} \VSState{(-2, 4)}{1}
\VSState{(-5, -4)}{5}\VSState{(5, -4)}{2}\VSState{(2, -3)}{3}\VSState{(-2, -3)}{4}
\EdgeR{p}{0}{pc} \EdgeR{1}{p}{cp}
\EdgeL{p}{2}{pa} \EdgeL{5}{p}{ap}
\EdgeR{p}{3}{pb} \EdgeR{4}{p}{bp}
\ZZEdge{0}{1}{}\ZZEdge{2}{5}{}\ZZEdge{3}{4}{}
\end{VCPicture}}
\caption{}\label{RauzyPicture}
\end{center}
\end{figure}
\end{enumerate}
\pfk
\begin{lmm}\label{Lemma3} Let $\mathbf{u}$ be an infinite word over a~ternary
alphabet with the complexity $\mathcal{C}(n) = 2n+1$ for
any $n \in \mathbb{N}$ and with the language $\mathcal{L}(\mathbf{u})$ closed under reversal.
Let $p \in \mathcal{L}(\mathbf{u})$ be a~palindrome with $\#{\rm Lext}(p)=2$ and let ${\mathcal P}(|p|)+{\mathcal P}(|p|+1)=4$.
\begin{enumerate}
\item  If $p$ has no palindromic extension, then $p$ is a~maximal
left special factor and there exists a~palindrome $q$ of the same
length such that $q$ has two palindromic extensions.
 \item If $p$ has two palindromic extensions, then there exists a~palindrome $q$ of
the same length such that $q$ has no palindromic extension and $q$
is a~maximal left special factor.
\end{enumerate}
\end{lmm}

\pf Denote ${\rm Lext}(p)=\{a,b\}$ and $|p|=n$. Since
$\Delta \mathcal{C}(n) = 2$, there exists a~factor $q\neq p$ of
the same length such that $\#{\rm Lext}(q)= 2$. Denote ${\rm
Lext}(q)=\{A,B\}$.
\begin{enumerate}
\item Assume that $p$ has no palindromic extension. The only factors
with length $n+2$ of the form $\ell_1p\ell_2$ where
$\ell_1,\ell_2 \in \mathcal{A} $ are $apb$ and $bpa$. The factor
$p$ is thus a~maximal left special factor. Let us recall that any prefix
of a~left special factor is again left special. This together with $\Delta
\mathcal{C}(n+1)=2$ implies that there exist two left special factors of
length $n+1$: $qx$ and $qy$ for $x\neq y$ with
\begin{equation}{\rm Lext}(qx)=\{A,B\}={\rm
Lext}(qy)\,.
\label{oba} \end{equation}
By reversal closeness of ${\mathcal L}(\mathbf u)$,
we obtain that ${\rm Lext}(\overline{q}A)=\{x,y\}$ and ${\rm
Lext}(\overline{q}B)=\{x,y\}$. Since there are no other left special factors besides $qx$ and $qy$ in $\mathcal{L}_{n+1}(\mathbf{u})$, necessarily $\overline{q}=q$ and
$\{A,B\}=\{x,y\}$. Because of \eqref{oba}, we deduce that
both $xqx$ and $yqy$ belong to the language $\mathcal{L}(\mathbf{u})$, i.e., the palindrome $q$ has two
palindromic extensions.

\item Suppose that $p$ has two palindromic extensions $apa$ and
$bpb$. In the Rauzy graph $\Gamma_n$, the  bispecial factor $p$
has the indegree and outdegree $2$, the left special factor $q$
has the indegree $2$ and the right special factor $\overline{q}$ has the outdegree
$2$. Moreover, the palindromes of length $n$ are exactly the vertices lying on the axis of symmetry and
the palindromes of length $n+1$ are exactly the edges crossing the axis. These facts together with ${\mathcal P}(n)+{\mathcal P}(n+1)=4$ imply that
the Rauzy graph $\Gamma_n$ can only look as depicted in Figure~\ref{RauzyPicture2}.
\begin{figure}[!h]
\begin{center}
\MediumPicture\VCDraw{%
\begin{VCPicture}{(-5,-5)(5,2)}
\ChgStateLineWidth{0.5}\StateVar[p]{(0, 0)}{p} \StateVar[q]{(-2, -4)}{q} \StateVar[\overline{q}]{(2, -4)}{q2}
\EdgeR{p}{q}{} \EdgeR{q2}{p}{}
\LoopN{p}{}{}\ArcR{q}{q2}{}\ArcR{q2}{q}{}
\end{VCPicture}}
\caption{}
\label{RauzyPicture2}
\end{center}
\end{figure}
Note that $q$ and $\overline{q}$ may coincide.
Let us first show, that necessarily  $apb \in
\mathcal{L}(\mathbf{u})$. If not, then it is impossible in
$\Gamma_n$ to leave the cycles in which only the vertex $p$ has
the indegree or the outdegree bigger that 1. It means that the word $\mathbf u$ is eventually periodic - a contradiction with the complexity.
Thus $apb, bpa  \in  \mathcal{L}(\mathbf{u})$. Consequently both $pa$ and $pb$
are left special factors of length $n+1$ with ${\rm
Lext}(pa)=\{a,b\} = {\rm Lext}(pb)$. Since $\Delta
\mathcal{C}(n+1) = 2$, no other left special factor of the same
length exists. Thus $q$ is the maximal left special and there
exists a~unique letter $x$ and a~unique letter $y$ such that $Aqx$
and $Bqy$ belong to the language ${\mathcal L}(\mathbf{u})$ and $x\neq y$. It
implies that ${\rm Lext}(\overline{q})=\{x,y\}$, i.e., the factor
$\overline{q}$ is a~left special factor of length~$n$, and therefore,
$\overline{q} = q$ and $\{A,B\}=\{x,y\}$. The Rauzy graph
$\Gamma_n$ has two vertices with indegrees $>1$ - the bispecial factors $p$ and $q$, see Figure~\ref{RauzyPicture3}.
\begin{figure}[!h]
\begin{center}
\MediumPicture\VCDraw{%
\begin{VCPicture}{(-5,-5)(5,2)}
\ChgStateLineWidth{0.5}\StateVar[p]{(0, 0)}{p} \StateVar[q]{(0, -4)}{q}
\ArcR{p}{q}{} \ArcR{q}{p}{}
\LoopN{p}{}{}\LoopS{q}{}{}
\end{VCPicture}}
\caption{}
\label{RauzyPicture3}
\end{center}
\end{figure}
Since $q$ is a~maximal left special factor, we have two disjoint
possibilities:  the first one is that $xqx$ and $yqy$ belong to
the language, the second one is that $xqy$ and $yqx$ belong to the
language. But the first possibility implies that in the Rauzy graph
$\Gamma_n$, it is impossible to leave the cycles containing
only one bispecial factor $q$ - a contradiction. Therefore the
second situation occurs and $q$ has no palindromic extension.
\end{enumerate}
\pfk

\bigskip

\pf[\it Proof of theorem \ref{kdyP}]
We will proceed by mathematical induction on $n$.
Obviously, ${\mathcal P}(0)=1$ and ${\mathcal P}(1)=3$. Assume that
${\mathcal P}(n)+{\mathcal P}(n+1)=4$ for some $n \geq 0$.
Let $p\in \mathcal{L}_{n}(\mathbf{u})$ be a~palindrome with
zero or two palindromic extensions. According to Lemma~\ref{Lemma3} there
exists a~palindrome $q$ of the same length, which is a~left special
factor as well. Since $\Delta \mathcal{C}(n)=2$, all other factors
including palindromes in $\mathcal{L}_{n}(\mathbf{u})$ have a~unique left extension.
According to Lemma \ref{Lemma1}, these
palindromes have a~unique palindromic extension. By Lemma~\ref{Lemma3}, the palindromes $p$ and $q$ together have
two palindromic extensions. Therefore, the number of palindromic
extensions of all palindromes in $\mathcal{L}_{n}(\mathbf{u})$
together is equal to the number of palindromes of length $n$.
Since every palindrome of length $n+2$ is a~palindromic extension of a~palindrome of length $n$,
we obtain ${\mathcal P}(n+1)+{\mathcal P}(n+2)=4$.
\pfk
\bigskip

\noindent{\it Proof of theorem \ref{kdyPE} } is a~direct consequence of the
previous lemmas.

\section{Counterexamples}\label{counterex}
In this last section, we will show that for ternary words, unlike binary words,
the properties ${\mathcal P}$ and ${\mathcal {PE}}$ are not equivalent and we will provide counterexamples to reversed implications in Theorems~\ref{kdyP} and~\ref{kdyPE}.

We have seen that for the computation of the first difference of complexity $\Delta {\mathcal C}(n)$, an important role is played by left and right special factors. See Formula~\eqref{eq:delta}.
In the sequel, it will be helpful to use a~formula for the second difference of complexity $\Delta^2{\mathcal C}(n)$, introduced by Cassaigne~\cite{Ca}. Let us explain that for the computation of $\Delta^2{\mathcal C}(n)$, bispecial factors are crucial.
Since every factor of length $n+2$ can be written as $xwy$, where
$x,y \in\mathcal{A}$ and $w \in \mathcal{L}(\mathbf u)$, it holds
$$\mathcal{C}(n+2) = \sum_{w\in \mathcal{L}_n(\mathbf{u})}\# \{ xwy\mid xwy \in
\mathcal{L}(\mathbf{u})\},$$ and similarly,
$$\mathcal{C}(n+1) = \sum_{w\in \mathcal{L}_n(\mathbf{u})}\# {\rm Lext}(w)  = \sum_{w\in
\mathcal{L}_n(\mathbf{u})}\# {\rm Rext}(w).$$ The second difference of complexity
$\Delta^2\mathcal{C}(n)=\Delta\mathcal{C}(n+1)-
\Delta\mathcal{C}(n) =\mathcal{C}(n+2) - 2\mathcal{C}(n+1)
+\mathcal{C}(n)$ may be obtained as follows
\begin{equation}\label{eq:delta2}
\Delta^2\mathcal{C}(n) =\sum_{w\in \mathcal{L}_n(\mathbf{u})}\Bigl(\# \{ xwy\mid xwy \in
\mathcal{L}(\mathbf{u})\} - \# {\rm Lext}(w)  - \# {\rm Rext}(w) +
1 \Bigr).
\end{equation}
Denote by $$b(w):= \# \{ xwy\mid xwy \in
\mathcal{L}(\mathbf{u})\} - \# {\rm Lext}(w)  - \# {\rm Rext}(w) +
1.$$ The number $b(w)$ is called the {\em bilateral order} of the factor
$w$. It is readily seen that if $w$ is not a~bispecial factor, then $b(w)=0$.
Bispecial factors will be distinguished according to their bilateral order in the following way
\begin{itemize}
\item if $b(w)>0$, then we call $w$ a~{\em strong} bispecial factor,
\item if $b(w)<0$, then we call $w$ a~{\em weak} bispecial factor,
\item if $b(w)=0$ and $w$ is bispecial, then we call it {\em ordinary}.
\end{itemize}
Evidently, for the value of $\Delta^2\mathcal{C}(n)$, only strong and weak bispecial factors are of importance.
\begin{rmrk}\label{useful_pal}
If $p$ is a~palindromic factor of a~reversal closed language ${\mathcal L}(\mathbf u)$, then $\#\{ xpy\mid xpy \in
\mathcal{L}(\mathbf{u})\}$ and the number of palindromic extensions of $p$ in $\mathbf u$ have the same parity. Moreover, $ \# {\rm Lext}(p)  = \# {\rm
Rext}(p)$. Therefore, the following simple observation holds
$$p \ \ \text {has a~unique palindromic extension in $\mathbf u$}\quad  \Longrightarrow \quad \text {
$b(p)$ is even}.$$
\end{rmrk}

\subsection{${\mathcal P}$ and ${\mathcal{PE}}$ are not equivalent}
The construction of a~ternary infinite word $\mathbf v$ with the desired properties is inspired by Arnoux and Rauzy \cite{ArRo} and Rote~\cite{Ro}.
Let $\mathbf v$ be the ternary infinite word defined by ${\mathbf v}=\Psi(\mathbf u)$, where
$\Psi:\{a,b\}^{*}\to \{0,1,2\}^{*}$ is the morphism given by
\begin{equation}\label{Psi}
\Psi(a)=12 \quad \text{and} \quad \Psi(b)=100,
\end{equation}
and $\mathbf u$ is the fixed point of the morphism $\varphi: \{a,b\}^{*}\to \{a,b\}^{*}$
defined by
\begin{equation}\label{varphi}
\varphi(a)=abbabba, \quad \varphi(b)=aba.
\end{equation}
In the sequel, we will show that $\mathbf v$ satisfies ${\mathcal P}$, but does not satisfy ${\mathcal {PE}}$. We will proceed in two steps. First, we will study several properties of the binary infinite word $\mathbf u$. Second, we will prove, using the properties of $\mathbf u$, that $\mathbf v$ satisfies $\mathcal P$, but does not satisfy ${\mathcal {PE}}$.
\subsubsection*{Step 1:} Let us show that the binary word $\mathbf u$ being the fixed point of the morphism $\varphi$ given  in~\eqref{varphi} has the language ${\mathcal L}({\mathbf u})$ closed under reversal and let us provide the list of all weak and strong bispecial factors of~$\mathbf u$.

Let us start with an important observation.
\begin{obs}\label{obs_decomp_varphi}
Every factor $v$ of $\mathbf u$ can be decomposed as $v=v^{(0)}v^{(1)}\dots v^{(m)}, \ m \geq 1,$ so that
$v^{(i)}\in \{aba, abbabba\}$ for $i \in \{1, \dots, m-1\}$, $v^{(0)}$ is a~proper suffix of $abbabba$ and $v^{(m)}$ is a~proper prefix of $abbabba$. Obviously, for every such decomposition, there exists $\tilde v \in \{a,b\}^{*}$ satisfying
\begin{equation}\label{decomp_varphi}
v=v^{(0)}\varphi(\tilde v)v^{(m)}.
\end{equation}
If the decomposition is unique, the corresponding $\tilde v$ is necessarily a~factor of $\mathbf u$.
\end{obs}
An essential role for the description of bispecial factors and palindromes in $\mathbf u$ is played by the map $T: \{a,b\}^{*} \to \{a,b\}^{*}$
defined by
\begin{equation}\label{T}
T(w)=ba\varphi(w)ab \quad \text{for every $w \in \{a,b\}^{*}$}.
\end{equation}
Let us summarize the properties of $T$ in the following lemma.
\begin{lmm}\label{Tproperties}
Let $T$ be the map defined in~\eqref{T}. Then, for every $w \in \{a,b\}^{*}$ and for all $c,d \in \{a,b\}$, it holds \begin{enumerate}[a)]
\item if $w$ is a~palindrome, then $T(w)$ is a~palindrome,
\item $cwd$ is a~factor of $\mathbf u$ if and only if $cT(w)d$ is a~factor of $\mathbf u$, in particular, if $w$ is a~factor of $\mathbf u$, then $T(w)$ is a~factor of $\mathbf u$.
\end{enumerate}
\end{lmm}
\begin{proof}
\begin{enumerate}[$a)$]
\item Since $\varphi(a)=aba$ and $\varphi(b)=abbabba$ are palindromes, it implies that $\varphi(w)$ is a~palindrome, thus $T(w)=ba\varphi(w)ab$ is a~palindrome, too.
\item $(\Rightarrow)$: If $awb$ is a~factor of $\mathbf u$, then $\varphi(awb)$ is in ${\mathcal L}(\mathbf u)$. As $aT(w)b$ is a~factor of $\varphi(awb)=aba\varphi(w)abbabba$, it follows that $aT(w)b$ is also a~factor of $\mathbf u$. The proofs for the other cases $awa, bwa, bwb$ are similar. $(\Leftarrow)$: Let $aT(w)b$ be a~factor of $\mathbf u$. It is readily seen that the unique decomposition of the form~\eqref{decomp_varphi} of $aT(w)b$ is $aT(w)b=\varphi(aw)abb$. Since $abb$ is a~prefix of $\varphi(b)$, but not of $\varphi(a)$, it follows that $awb \in {\mathcal L}(\mathbf u)$. The proofs for the other cases $aT(w)a, bT(w)a, bT(w)b$ are analogous.
\end{enumerate}
\end{proof}
\begin{rmrk} Lemma~\ref{Tproperties} has several useful consequences.
\begin{enumerate}
\item According to Lemma~\ref{Tproperties}, the language ${\mathcal L}(\mathbf u)$ contains infinitely many palindromes. Together with the primitivity of the substitution $\varphi$, thus the uniform recurrence of $\mathbf u$, it implies that the language ${\mathcal L}(\mathbf u)$ is closed under reversal.
\item For any factor $w \in {\mathcal L}(\mathbf u)$, its bilateral order $b(w)=b(T(w))$ by Item~$b)$ of Lemma~\ref{Tproperties}.
\item If $w$ is a~palindrome in ${\mathcal L}(\mathbf u)$, then $T(w)$ is a~palindrome with the same number of palindromic extensions by Lemma~\ref{Tproperties}.
\end{enumerate}
\end{rmrk}
Since the word $\mathbf u$ is built from the factors $abbabba$ and $aba$, it is clear that the words $$aaa, bbb, abab, baba, aabbaa, babbab$$ are not in ${\mathcal L}(\mathbf u)$. Observing then the prefix of $\mathbf u$
$${\mathbf u}=abbabbaabaabaabbabbaabaaba\dots,$$
it follows that the only left special factors of length $\leq 4$ are: $\varepsilon, a,b,ab,ba,abb,baa,abba,baab$; among them, only $\varepsilon$ and $baab$ are strong bispecial factors and only $abba$ is a~weak bispecial factor.
\begin{lmm}\label{bs_u}
For every bispecial factor $v\in {\mathcal L}(\mathbf u)$ of length at least $5$, there exists a~factor $w \in {\mathcal L}(\mathbf u)$ such that $v=T(w)$.
Moreover, $b(w)=b(T(w))$.
\end{lmm}
\begin{proof}
Every prefix of a~left special factor is left special, too.
Since $abba$ and $baab$ are the only left special factors of length $4$ and $abba$ is a~weak bispecial factor, thus cannot be extended to the right staying left special, we learn that every bispecial (thus left special) factor $v$ of length $\geq 5$ has to start in $baab$. Since the language ${\mathcal L}(\mathbf u)$ is closed under reversal, the bispecial (thus right special) factor $v$ has to end in $baab$. Then, it is clear from the form of the substitution $\varphi$ that $v=ba\varphi(w)ab$ is a~unique decomposition of the form~\eqref{decomp_varphi} of $v$.  Thus, by Observation~\ref{obs_decomp_varphi}, $w$ is a~factor of $\mathbf u$ such that $v=T(w)$.
The last statement is a~consequence of Item~$b)$ of Lemma~\ref{Tproperties}.
\end{proof}
As a~consequence of Lemma~\ref{bs_u}, we obtain the set of all strong bispecial factors
\begin{equation}\label{strong}\{V^{(n)}\mid n \in \mathbb N\}, \quad \text{where $V^{(0)}=\varepsilon$ and $V^{(n)}=T(V^{(n-1)})$ for $n \geq 1$},
\end{equation} and the set of all weak bispecial factors
\begin{equation}\label{weak}\{U^{(n)}\mid n \in \mathbb N, \ n \geq 1\}, \quad \text{where $U^{(1)}=abba$ and $U^{(n)}=T(U^{(n-1)})$ for $n \geq 2$}.\end{equation}
It is easy to see that $b(V^{(0)})=b(\varepsilon)=1$ and $b(U^{(1)})=b(abba)=-1$.
Item $b)$ of Lemma~\ref{Tproperties} implies that $b(V^{(n)})=1$ and $b(U^{(n)})=-1$ for all $n$.
Moreover, by Item~$a)$ of Lemma~\ref{Tproperties}, they are palindromes.

\subsubsection*{Step 2:}
We may now study the ternary word $\mathbf v=\Psi(\mathbf u)$ defined in~\eqref{Psi}.
In the sequel, it will be shown that
\begin{enumerate}
\item the language ${\mathcal L}(\mathbf v)$ is closed under reversal,
\item the complexity of $\mathbf v$ is ${\mathcal C}(n)=2n+1$ for all $n \in \mathbb N$,
\item the language ${\mathcal L}(\mathbf v)$ contains infinitely many distinct palindromes that do not have a~unique palindromic extension.
\end{enumerate}
When proven, the statements $(1)$ and $(2)$ imply that {\bf the property $\mathcal P$ holds} (by Theorem~\ref{kdyP}) and the statement $(3)$ has as a~consequence that {\bf the property ${\mathcal {PE}}$ does not hold}.\\

\vspace{0.1cm}
\noindent{\it Proof of Step 2:}

Let us start with a~similar observation as Observation~\ref{obs_decomp_varphi}.
\begin{obs}\label{obs_decomp_Psi}
Every factor $v$ of $\mathbf v$ can be decomposed as $v=v^{(0)}v^{(1)}\dots v^{(m)}, \ m \geq 1,$ so that
$v^{(i)}\in \{12, 100\}$ for $i \in \{1, \dots, m-1\}$, $v^{(0)}$ is a~proper suffix either of $12$ or of $100$ and $v^{(m)}$ is a~proper prefix of $100$. Obviously, for every such decomposition, there exists $\tilde v \in \{a,b\}^{*}$ such that
\begin{equation}\label{decomp_Psi}
v=v^{(0)}\Psi(\tilde v)v^{(m)}.
\end{equation}
If the decomposition is unique, the corresponding $\tilde v$ is necessarily a~factor of $\mathbf u$.
\end{obs}
The crucial tool for the proof of $(1)-(3)$ is the map $H: \{a,b\}^{*} \to \{0,1,2\}^{*}$ defined by
\begin{equation}\label{hatT} H(w)=\Psi(w)1 \quad \text{for every $w \in \{a,b\}^{*}$}.\end{equation}
Its properties are stated in the following lemma.
\begin{lmm}\label{hatTproperties}
Let $H$ be the map defined in~\eqref{hatT}. Then it holds for every $w \in \{a,b\}^{*}$
\begin{enumerate}[a)]
\item if $w$ is a~factor of $\mathbf u$, then $H(w)$ is a~factor of $\mathbf v$,
\item if $w$ is a~palindrome, then $H(w)$ is a~palindrome,
\item if $w$ is a~factor of $\mathbf u$, then $b(w)=b(H(w))$.
\end{enumerate}
\end{lmm}
\begin{proof}
\begin{enumerate}[$a)$]
\item There exists a~letter $x \in \{a,b\}$ such that $wx \in {\mathcal L}(\mathbf u)$. Then $\Psi(wx)$ is a~factor of $\mathbf v=\Psi(\mathbf u)$ and $\Psi(wx)$ contains $H(w)=\Psi(w)1$.
\item It suffices to notice that $1^{-1}\Psi(a)1=\overline{\Psi(a)}$ and $1^{-1}\Psi(b)1=\overline{\Psi(b)}$, where $1^{-1}\Psi(a)1$ is the word obtained when the prefix $1$ is cut from $\Psi(a)1$.
\item The statement will be proven if we show that the relation between the extensions of $w$ and $H(w)$ is as follows
$$awa \in {\mathcal L}(\mathbf u) \quad \Leftrightarrow \quad 2H(w)2 \in {\mathcal L}(\mathbf v)$$
$$awb \in {\mathcal L}(\mathbf u) \quad \Leftrightarrow \quad 2H(w)0 \in {\mathcal L}(\mathbf v)$$
$$bwa \in {\mathcal L}(\mathbf u) \quad \Leftrightarrow \quad 0H(w)2 \in {\mathcal L}(\mathbf v)$$
$$bwb \in {\mathcal L}(\mathbf u) \quad \Leftrightarrow \quad 0H(w)0 \in {\mathcal L}(\mathbf v)$$
$(\Rightarrow)$: If $awb \in {\mathcal L}(\mathbf u)$, then $\Psi(awb)=12\Psi(w)100$ is a~factor of $\mathbf v$ and $\Psi(awb)$ contains $2H(w)0$. The proofs for the other cases $awa, bwa, bwb$ are similar. $(\Leftarrow)$: It is easy to see that $2H(w)0=2\Psi(w)10$ is a~unique decomposition of $2H(w)0$ of the form~\eqref{decomp_Psi}. Moreover, since $2$ is a~suffix of $\Psi(a)$, but not of $\Psi(b)$, and $10$ is a~prefix of $\Psi(b)$, but not of $\Psi(a)$, it follows that $awb$ is a~factor of $\mathbf u$. The proofs for the other cases $2H(w)2, 0H(w)2, 0H(w)0$ are analogous.
\end{enumerate}
\end{proof}
\begin{enumerate}
\item
According to its construction, the word $\mathbf v$ is uniformly recurrent. Using Items~$a)$ and~$b)$ of Lemma~\ref{hatTproperties}, it is clear that ${\mathcal L}(\mathbf v)$ contains infinitely many distinct palindromes. Relating these two facts, ${\mathcal L}(\mathbf v)$ is closed under reversal.
\item In order to describe all strong and weak bispecial factors, the following lemma is helpful.
However, it is useful to notice first that the only left special factors of length $\leq 2$ are: $\varepsilon, 0, 1, 10, 12$. Among them, the only strong bispecial factor is $1$ and the only weak bispecial factor is $0$.
\begin{lmm}\label{bs_v}
Let $v$ be a~bispecial factor of $\mathbf v$ of length $\geq 3$. There exists a~factor $w$ of $\mathbf u$ such that
$v=H(w)$. Moreover, $b(w)=b(H(w))$.
\end{lmm}
\begin{proof}
Since every prefix of a~bispecial factor is left special, $v$ has to start in $1$. Since the language ${\mathcal L}(\mathbf v)$ is closed under reversal and $v$ is right special, $v$ has to end in $1$. Then, observing the morphism $\Psi$, $v=\Psi(w)1$ is a~unique decomposition of $v$ the form~\eqref{decomp_Psi}. Thus, by Observation~\ref{obs_decomp_Psi}, $w$ is a~factor of $\mathbf u$ satisfying $v=H(w)$.
The last statement follows by Item~$c)$ of Lemma~\ref{hatTproperties}.
\end{proof}
By Lemma~\ref{bs_v} and since $0$ is the only weak and $1$ the only strong bispecial factor of length $\leq 2$, we obtain the set of all strong bispecial factors of $\mathbf v$ (recall that $V^{(n)}$ and $U^{(n)}$ are defined in~\eqref{strong} and in~\eqref{weak}, respectively)
$$\{\hat V^{(n)}\mid n \in \mathbb N\}, \quad \text{where $\hat V^{(n)}=H(V^{(n)})$},$$
and the set of all weak bispecial factors of $\mathbf v$
$$\{\hat U^{(n)}\mid n \in \mathbb N\}, \quad \text{where $\hat U^{(0)}=0$ and $\hat U^{(n)}=H(U^{(n)})$ for $n \geq 1$}.$$
Since the factors $U^{(1)}=abba$ and $V^{(1)}=baab$ consist of the same ``hand'' of letters, it follows by the definition of $V^{(n)}$ and $U^{(n)}$ that $|V^{(n)}|_a=|U^{(n)}|_a$ and $|V^{(n)}|_b=|U^{(n)}|_b$, where $|w|_a$ denotes the number of letters $a$ occurring in a~word $w$. Therefore, we deduce that $|\hat V^{(n)}|=|\hat U^{(n)}|$ and by Lemma~\ref{bs_v}, it holds $b(\hat V^{(n)})+b(\hat U^{(n)})=b(V^{(n)})+b(U^{(n)})=0$.
By Formula~\eqref{eq:delta2}, we have $\Delta^2 {\mathcal C}(n)\equiv 0$, and since ${\mathcal C}(0)=1$ and ${\mathcal C}(1)=3$, it follows that $C(n)=2n+1$ for every $n \in \mathbb N$.
\item
The strong bispecial factors $\hat V^{(n)}$ are palindromes by Item~$b)$ of Lemma~\ref{hatTproperties}.
Since $b(\hat V^{(1)})=b(1)=1$, we deduce using Item~$c)$ of Lemma~\ref{hatTproperties} that $b(\hat V^{(n)})=1$ for all $n \in \mathbb N$. Applying Remark~\ref{useful_pal}, the palindromes $\hat V^{(n)}$ do not have a~unique palindromic extension. Using similar arguments, the palindromes $\hat U^{(n)}$ do not have a~unique palindromic extension either.
\end{enumerate}
%%%%%%%%%%%%%%%%%%%%%%%%%%%%%%%%%%%%%%%%%%%%%%%%%%%%%%%%%%%%%%%%%%%%%%%%%%%%%%%%%%%%%%%%%%%%%%%%%%%%%%%%%%%%%%%%%%%%%%
\subsection{Implications in Theorems~\ref{kdyP} and~\ref{kdyPE} are irreversible}
In order to show that the implications in Theorems~\ref{kdyP} and~\ref{kdyPE} are irreversible, we will construct an infinite ternary word $\mathbf U$ whose language ${\mathcal L}(\mathbf U)$ is closed under reversal and such that on one hand, $\mathbf U$ has the property ${\mathcal {PE}}$, consequently $\mathbf U$ has the property $\mathcal P$, too, on the other hand, the complexity ${\mathcal C}(n)$ of $\mathbf U$ does not satisfy ${\mathcal C}(n)=2n+1$ for all $n \in \mathbb N$.

Denote by $\mathbf U$ the infinite ternary word being the fixed point of the morphism $\Phi: \{A,B,C\}^{*} \to \{A,B,C\}^{*}$ defined by
\begin{equation}\label{U}
\Phi(A)=ABA, \quad \Phi(B)=CAC, \quad \Phi(C)=ACA.
\end{equation}
We will not provide a~detailed proof of the announced properties, but only a~helpful hint for the reader.
Observing the substitution $\Phi$, it is obvious that the image of a~palindrome is again a~palindrome. Therefore, ${\mathcal L}(\mathbf U)$
contains infinitely many palindromes. Together with the uniform recurrence of $\mathbf U$, it implies that the language ${\mathcal L}(\mathbf U)$ is closed under reversal.
In addition, every palindrome $p$ is a~central factor of $\Phi^2(p)$, i.e., there exists $w \in \{A,B,C\}^{*}$ such that $\Phi^2(p)=wp\overline{w}$. In particular,
$(\Phi^{2n}(A))$ is a~sequence of palindromes with $A$ as a~central factor, $(\Phi^{2n}(B))$ is a~sequence of palindromes with $B$ as a~central factor,
$(\Phi^{2n}(C))$ is a~sequence of palindromes with $C$ as a~central factor, and $(\Phi^{2n}(AA))$ is a~sequence of palindromes of even length.
It is easy to see that every palindrome is a~central factor of one of the above families, thus the property ${\mathcal {PE}}$ holds.

Concerning the complexity, we have $${\mathcal L}_3(\mathbf U)=\{AAB, BAA, AAC, CAA, ABA, ACA, CAC, BAC, CAB\},$$ hence ${\mathcal C}(1)=3, \ {\mathcal C}(2)=5, \ {\mathcal C}(3)=9$.
Thus, it does not hold ${\mathcal C}(3)=2\cdot 3+1$. In fact, $\Delta {\mathcal C}(n)\not = 2$ for infinitely many $n \in \mathbb N$.
%%%%%%%%%%%%%%%%%%%%%%%%%%%%%%%%%%%%%%%%%%%%%%%%%%%%%%%%%%%%%%%%%%%%%%%%%%%%%%%%%%%%%%%
%%%%%%%%%%%%%%%%%%%%%%%%%%%%%%%%%%%%%%%%%%%%%%%%%%%%%%%%%%%%%%%%%%%%%%%%%%%%%%%%%%%%%%%%%%
\section*{Acknowledgements}
The authors acknowledge financial support by the grants MSM6840770039 and LC06002 of the
Ministry of Education, Youth, and Sports of the Czech Republic.
%%%%%%%%%%%%%%%%%%%%%%%%%%%%%%%%%%%%%%%%%%%%%%%%%%%%%%%%%%%%%%%%%%%%%%%%%%%%%%%%%%%%%%%%%%%%%%%%%%%%%%%%%%%%%%%%%%%%%%%%


\begin{thebibliography}{10}
\bibitem{ArRo} P. Arnoux, G. Rauzy, Repr\'esentation g\'eom\'etrique de suites de complexit\'e $2n+1$. Bull. Soc. Math. France {\bf 119} (1991) 199-215.
\bibitem{BaMaPe} P. Bal\'a\v zi, Z. Mas\'akov\'a,
E. Pelantov\'a, Factor versus palindromic complexity of
uniformly recurrent infinite words. Theoret. Comput. Sci. {\bf 380} (2007) 266-275.
\bibitem{BaPeSt} L{\!'}. Balkov\'a, E. Pelantov\'a, W. Steiner, Sequences with Constant Number
of Return Words. Monatshefte f\"{u}r Mathematik, {\bf 155(3-4)} (2008) 251-263.
\bibitem{Ca} J. Cassaigne, Complexity and special factors.
Bull. Belg. Math. Soc. Simon Stevin 4 {\bf 1} (1997) 67-88.
\bibitem{DrPi} X. Droubay, G. Pirillo, Palindromes and Sturmian words. Theoret.
Comput. Sci. {\bf 223} (1999) 73-85.
\bibitem{GlJuWiZa} A. Glen, J. Justin, S. Widmer, L. Q. Zamboni, Palindromic richness. Eur. J. Comb. {\bf 30} (2009) 510-531.
\bibitem{HeMo} G. A. Hedlund, M. Morse, Symbolic dynamics II - Sturmian trajectories. Amer. J. Math.  {\bf 62} (1940), 1-42.
\bibitem{Ro} G. Rote, Sequences with subword complexity $2n$. Journal of Number Theory {\bf 46} (1993) 196-213.
\bibitem{Vu} L. Vuillon, A~characterization of Sturmian words by
return words. Eur. J. Comb. {\bf 22} (2001)  263-275.




\end{thebibliography}
\end{document}